\documentclass[11pt]{article}
\usepackage{graphicx,epsfig}

\begin{document}

\begin{center}
{\LARGE  On  stabilizability and exact observability of stochastic
systems with their applications \footnote{This paper was completed
in 2001  and submitted to Automatica as a regular paper for possible
publication, which has been  published as a brief paper in
Automatica, 40: pp.87-94, 2004. This material will help interesting
reader read ``Automatica,  40, 87-94, 2004" with ease. }}

\vspace{0.5cm}

{Weihai Zhang }${}^1$
 and { Bor-Sen Chen}${}^2$

 ${}^1$ {College of Electrical Engineering and Automation,}

 {Shandong
University of Science and Technology,}

 {Qingdao 266590, China}

{Email: w\_hzhang@163.com}

${}^2$ {Department of Electrical Engineering, National Tsing Hua
University,}

{Hsinchu 30013, Taiwan}

{Email: bschen@ee.nthu.edu.tw}
\end{center}

\begin{center}
{\large Abstract}
\end{center}

This paper discusses  the stabilizability, weak stabilizability,
exact observability and robust quadratic stabilizability  of linear
stochastic  control systems. By means  of the spectrum technique  of
the generalized Lyapunov operator, a necessary and sufficient
condition is given for stabilizability  and  weak stabilizability of
stochastic systems, respectively.
 Some new concepts
  called  unremovable spectrums,
strong solutions, and weakly feedback  stabilizing solutions are
introduced.  An unremovable spectrum  theorem is given, which
generalizes the corresponding theorem of deterministic systems to
stochastic systems. A stochastic Popov-Belevith-Hautus (PBH)
criterion for exact observability is obtained. For  applications, we
give a comparison theorem for generalized  algebraic Riccati
equations (GAREs), and two results on Lyapunov-type equations are
obtained, which improve the previous works. Finally, we also discuss
robust quadratic stabilization of uncertain stochastic systems, and
a necessary and sufficient condition is given for quadratic
stabilization via a linear matrix inequality (LMI).


\section{Introduction}

Stabilizability and observability is essential and important
concepts in modern control theory, especially, in system analysis
and  synthesis. Stochastic  stabilizability (in mean square sense)
plays a crucial role and  is an essential assumption in many
problems, such as  infinite horizon stochastic optimal control
problem (J. L. Willems \& J.~C.~Willems,~1976; Ichikawa, 1979; Ait
Rami \& Zhou,~2000), robust and stochastic $H^\infty$ problems
(Ugrinovskii, 1998; Hinrichsen \& Pritchard,~1998), filtering
problems (Bensoussan, 1992, and the reference therein), control and
stabilization problems for jump system ( Ghaoui \& Ait Rami,~1998,
Gao et al.,~2001 and the references therein). In this paper, we
mainly study the stabilizability and exact observability of the
following linear   stochastic controlled system
\begin{equation}
\left\{
\begin{array}{l}
dx(t)=(Ax(t)+Bu(t))\,dt+(Cx+Du)\,dw(t), x(0)=x_0\in {\cal R}^n, \\
y=Qx,
\end{array}
\right.
\end{equation}
where $A,B,C,D, Q\in {\cal R}^{n\times n}\times {\cal R}^{n\times m}
\times{\cal R}^{n\times n}\times{\cal R}^{n\times m}\times{\cal
R}^{l\times n}$ are real constant matrices; $w(\cdot)$ is a standard
Wiener process with $w(0)=0$  defined on the filtered probability
space $(\Omega,{\cal F},{\cal P};{\cal F}_t)$ with ${\cal
F}_t=\sigma\{w(s),0\le s\le t\}$. Without loss of generality, we
assume $w(\cdot)$ is one-dimensional for simplicity.

On the stabilizability of system (1) (For short, we also call
$(A,B;C,D)$  stabilizable in the context), some results were
obtained, for example, Zhang (2000) proved that $(A,B;C,D)$ is
stabilizable  if and only if (iff)  the following   GARE
\begin{equation}
PA+A'P+C'PC+I-(PB+C'PD)(I+D'PD)^{-1}(B'P+D'PC)=0
\end{equation}
has a positive definite solution $P>0$.  By means of linear matrix
inequalities (LMIs) and Lyapunov-type inequalities, Ait Rami and
Zhou (2000) gave some other  necessary and sufficient conditions. In
addition, J. L. Willems and J. C. Willems (1976) also presented many
criteria to test the stabilizability for some special stochastic
systems with arbitrary noise intensities.

Exact observability was first introduced by Zhang (1998) by means of
finite-dimensional invariance subspace, and then by Liu (1999) from
the physical viewpoint, i.e., the zero output must correspond to
zero state,  which is the generalized version of complete
observability of deterministic linear system theory, and has many
applications in studying   GAREs, and stochastic stability; see
Zhang (1998) and Liu (1999).

In this present paper, we try to give a spectral  or eigenvalue
descriptions for the stabilizability and exact observability of
stochastic system (1). In order to illustrate our goal, let  us
recall some results with  deterministic controlled system
\begin{equation}
\left\{
\begin{array}{l}
{\dot x}=Ax+Bu, \ \ \ x(t_0)=x_0 , \nonumber \\
y(t)=Qx(t). \label {eq 1.3}
\end{array}
\right.
\end{equation}
It is well known that for complete controllability, complete
observability, stabilizability and detectability of system  (\ref{eq
1.3}), we have the following criteria:

1)\  $(A,B)$ is completely controllable  iff
 there does not exist a nonzero vector $\xi$ and $\lambda\in {\cal C}$   satisfying
\begin{equation}
{\xi}'A=\lambda {\xi}',\ \  {\xi}'B=0.
\label {eq quo}
\end{equation}

2) \ $(A,B)$ is stabilizable iff  there does not exist a nonzero
vector $\xi$, $Re(\lambda)\ge 0$ satisfying (\ref{eq quo}).

3)\  $(Q,A)$ is completely observable iff
 there does not exist a nonzero vector $\xi$ and $\lambda\in {\cal C}$   satisfying
\begin{equation}
A{\xi}=\lambda {\xi},\ \  Q{\xi}=0.
\label {eq qub}
\end{equation}

4)\ $(Q,A)$ is detectable  iff  there does not exist a nonzero
vector $\xi$, $Re(\lambda)\ge 0$ satisfying (\ref{eq qub}).

The above criteria for complete controllability(stabilizability, complete observability, detectability, resp.)   are  the so-called   Popov-Belevith-Hautus (PBH) Criteria,
and the $\lambda$ satisfying (4) can be  called an unremovable
spectrum of system (3); see Proposition 1 of this paper. So PBH
criterion tells us that system (3) is stabilizable iff its all
unremovable spectra  belong to left hand side complex plane. PBH
criterion is important in the pole assignment of linear systems.

In this paper, we try to develop an analogous theory in stochastic
systems. By defining the closed-loop operator ${\cal L}_K$,
 the necessary and sufficient conditions for stabilizability
and weak stabilizability of system (1) are presented via the
spectrum of operator ${\cal L}_K$. The proposed spectrum-based
technique makes it possible  to  generalize some concepts of
deterministic systems to  stochastic case.  For instance, we can
define weak stabilizability of stochastic systems, strong solution
of GARE, etc.. We obtain a stochastic unremovable spectrum theorem,
however, PBH criterion  is only a necessary but not a sufficient
condition for the stochastic stabilizability. It is very interesting
that PBH criterion  still holds for exact observability. Moreover,
we find that even if system (1) is exactly-terminal controllable
(Peng 1994), the spectrum  of system (1)  cannot be assigned
arbitrarily.  All these reveal the essential differences between the
deterministic and stochastic systems.

As our theoretical applications, some results on GAREs  and
Lyapunov-type equations are improved. This paper also discusses
the quadratic stabilization of stochastic uncertain systems, for which, a
necessary and sufficient condition is given.

This paper is organized as follows: In section 2, we discuss the
relation between stabilizability and spectrum. Some definitions are
introduced, then necessary and sufficient theorems are obtained for
the stabilizability and unremovable spectrum. The stochastic PBH
Criteria  are also presented for some special systems. In
section 3, a stochastic PBH  criterion is obtained for exact
observability, by duality, we also define stochastic detectability,
various implication  relations with complete observability, exact
observability, detectability and stochastic detectability have been
clarified.  In section 4, we deal with the weak stabilizability of
stochastic systems, and a necessary and sufficient theorem is given in
form of spectrum, and some sufficient conditions are given for the
weak stabilizability in LMIs and Lyapunov-type inequalities. In
section 5, some theoretical applications to GAREs  and Lyapunov-type
equations are developed.  Section 6 gives a necessary and sufficient
condition for quadratic stabilization of stochastic uncertain
systems. Section 7 concludes this paper with some comments.

\section{ Stabilizability and spectrum of stochastic systems }

For convenience, we adopt the following notations:

${\cal S}_n$: the set of all $n\times n$ symmetric matrices, its
entries may be complex numbers;

$A'$ ( $Ker(A)$): the  transpose ( kernel space) of the
 matrix $A$;

$A\ge 0 (A>0)$:  $A$  is a positive semidefinite (positive definite)
symmetric matrix;

$I$: identity matrix;

$\sigma(L)$: spectral set of the operator or matrix $L$;

${\cal C}^{-}({\cal C}^{-,0})$: the open left (closed left ) hand side complex plane.

In this section, we mainly study the stabilizability and spectral
properties of stochastic system (1). Firstly,  we  give the
following definitions:

{\bf Definition 1.} Stochastic system (1) is called stabilizable (in the mean square sense),
if there exists a feedback
control $u(t)=Kx(t)$, such that for any $x_0\in {\cal R}^n$, the closed-loop system
\begin{equation}
dx(t)=(A+BK)x(t)\,dt+(C+DK)x(t)\,dw(t), x(0)=x_0
\label {eq 78}
\end{equation}
is asymptotically mean square stable, i.e., $\lim_{t\to \infty}
E[x(t)x'(t)]=0$, where $K\in {\cal R}^{m\times n}$ is a constant
matrix.

Below, we define the spectrum of the closed-loop system (\ref{eq
78}).

{\bf Definition 2.} For any given feedback gain matrix $K$,  let
${\cal L}_K$ be a linear operator from ${\cal S}_n$ to ${\cal S}_n$
defined as follows:
$$
{\cal L}_K: X\in {\cal S}_n \longmapsto
(A+BK)X+X(A+BK)'+(C+DK)X(C+DK)'.
$$
The spectrum of ${\cal L}_K$ is the set defined by $\sigma ({\cal
L}_K)=\{\lambda\in {\cal C}:{\cal L}_K (X)=\lambda X, X\in {\cal
S}^n, X\ne 0\}$, which  is also called the spectrum of system
(\ref{eq 78}).

Now, we give the first theorem for stabilizability of system (1).

{\bf Theorem 1.} System (1) is stabilizable iff  there exists a
$K\in {\cal R}^{m\times n}$, such that the spectrum of (\ref{eq 78})
belongs to ${\cal C}^{-}$.

{\bf Proof.}   By definition, we only need to prove there exists a
$K\in {\cal R}^{m\times n}$, such that (\ref{eq 78}) is
asymptotically  mean square stable. Let $X(t)=E[x(t)x'(t)]$, where
$x(t)$ is the trajectory of (\ref{eq 78}), then by Ito's formula,
\begin{equation}
\left\{
\begin{array}{l}
\dot X(t)=(A+BK)X(t)+X(A+BK)'+(C+DK)X(C+DK)'={\cal L}_K(X(t)), \nonumber \\
X(0)=X_0=x_0x'_0. \label {eq 05}
\end{array}
\right.
\end{equation}
Since $X(t)\in {\cal S}_n$,  (\ref{eq 05}) is a symmetric
matrix-valued equation including $n(n+1)/2$ different variables. If
we write $X=(Ex_ix_j)_{n\times n}=(X_{ij})_{n\times n}$, and define
a map $\tilde{\cal L}$ from ${\cal S}_n$ to ${\cal R}^{\frac
{n(n+1)} 2}$ as follows:
$$
\tilde{X}=\tilde{\cal L}(X)=(X_{11},X_{12},\cdots, X_{1n}, X_{22},X_{23},\cdots, X_{2n},\cdots,
X_{n-1,n-1},X_{n-1,n},X_{nn})',
$$
then there exists a unique matrix\footnote {In order to give the
exact expression of $L(0)$, we developed a new technique called
${\mathcal H}$-representation technique in 2012; see W. Zhang and B.
S. Chen, ${\mathcal H}$-representation and applications  to
generalized Lyapunov equations and linear stochastic systems, IEEE
Trans. Automatic Control, 57(12): pp. 3009-3022, 2012.} $L(K)\in
{\cal R}^{\frac {n(n+1)} 2\times \frac {n(n+1)} 2}$, such that
equation (\ref{eq 05} ) is equivalent to
\begin{equation}
\dot{\tilde{X}}=L(K)\tilde{X},\ \ \  \tilde{X}(0)=\tilde {X_0}.
\label {eq 06}
\end{equation}
Obviously,
\begin{equation}
\lim_{t\to\infty} Ex(t)x'(t)=0\Longleftrightarrow\lim_{t\to\infty}\tilde{X}(t)=0
\Longleftrightarrow \sigma (L(K))\subset {\cal C}^- .
\end{equation}
By Definition 2, it is not difficult to prove  $\sigma (L(K))=\sigma
({\cal L}_K)$.  Hence,  from the above, the proof of Theorem 1 is
completed.

Below, we say $L(K)$ is a matrix induced by  ${\cal L}_K$.

{\bf Remark 1.} It is easily seen that the following operator
$$
{\cal L}^* _K: X\in {\cal S}_n\longmapsto X(A+BK)+(A+BK)'X+(C+DK)'X(C+DK)
$$
is the adjoint  operator of ${\cal L}_K$ with the inner product
$<X,Y>=trace(X^*Y)$ for any $X,Y\in {\cal S}_n$. So system (1) is
stabilizable iff there exists an $K\in {\cal R}^{m\times n}$, such
that $\sigma ({\cal L}^* _K) \subset {\cal C}^{-}$. As we limit the
coefficient matrices to be real, so $\sigma ({\cal L}^* _K)=\sigma
({\cal L}_K)$. If we denote the induced matrix of the operator
${\cal L}^* _K$ by $L^*(K)$, then from the proof of Theorem 1, we
have
$$
\sigma ({\cal L}^* _K)=\sigma ({\cal L}_K)=\sigma (L(K))=\sigma (L^*(K)),
$$
where any one of them can characterize the stabilizability of system (1).

To illustrate the meaning of the above notations, we give an example as follows:

{\bf Example 1:} In system (\ref{eq 78}), take $K=I$,
$$
A=\left[
\begin{array}{cc}
-3 &  1/2\\
-1 &  -1
\end{array}
\right],  \quad
B=\left[
\begin{array}{cc}
1 & 0\\
0 &  1
\end{array}
\right], \quad
C=\left[
\begin{array}{cc}
2 & 0\\
0 &  0
\end{array}
\right], \quad
D=\left[
\begin{array}{cc}
2 & 0\\
0 &  0
\end{array}
\right].
$$
It is easily tested that (\ref{eq 05}) becomes
$$
\left[
\begin{array}{cc}
\dot X_{11} &  \dot X_{12}\\
 \dot X_{12} &   \dot X_{22}
\end{array}
\right]=\left[
\begin{array}{cc}
 X_{12} &   -X_{11}+\frac {1} {2} X_{22}\\
 -X_{11}+\frac {1} {2} X_{22} &   -2X_{12}+X_{22}
\end{array}
\right],
$$
which is equivalent to
$$
\dot{\tilde{X}}=
\left[
\begin{array}{c}
\dot X_{11} \\
 \dot X_{12} \\
 \dot X_{22}
\end{array}
\right]=
L(K)\tilde{X}=
\left[
\begin{array}{ccc}
0 & 1 & 0\\
-1 &   0 & \frac {1} 2\\
0 & -2 & 1
\end{array}
\right]
\left[
\begin{array}{c}
 X_{11} \\
  X_{12} \\
  X_{22}
\end{array}
\right].
$$
So
$$
L(K)=\left[
\begin{array}{ccc}
0 & 1 & 0\\
-1 &   0 & \frac {1} 2\\
0 & -2 & 1
\end{array}
\right].
$$
Since ${\cal L}_{K=I}(X)=\lambda X$, $X\in {\cal S}^n$,  is equivalent to
$$
\left[
\begin{array}{ccc}
0 & 1 & 0\\
-1 &   0 & \frac {1} 2\\
0 & -2 & 1
\end{array}
\right]
\left[
\begin{array}{c}
 X_{11} \\
  X_{12} \\
  X_{22}
\end{array}
\right]
=\lambda
\left[
\begin{array}{c}
 X_{11} \\
  X_{12} \\
  X_{22}
\end{array}
\right],
$$
$\sigma ({\cal L}_{K=I})=\sigma (L(K=I))$. From ${\cal
L}^*_{K=I}(X)=\lambda X$, we have
$$
\left[
\begin{array}{ccc}
0 & -2 & 0\\
\frac {1}{2} &   0 & -1\\
0 & 1 & 1
\end{array}
\right]
\left[
\begin{array}{c}
 X_{11} \\
  X_{12} \\
  X_{22}
\end{array}
\right]
=\lambda
\left[
\begin{array}{c}
 X_{11} \\
  X_{12} \\
  X_{22}
\end{array}
\right].
$$
 Via a
simple computation, we have
$\sigma ({\cal K}^*_{K=I})=\sigma ({L}^*(K=I))=\sigma ({\cal L}_{K=I})=\sigma (L(K=I))
=\{\lambda _1,\lambda _2,\lambda _3\}$, where
$$
L^*(K=I)=
\left[
\begin{array}{ccc}
0 & -2 & 0\\
\frac {1}{2} &   0 & -1\\
0 & 1 & 1
\end{array}
\right],
$$
and $\lambda _i$ are the roots of the following characteristic
polynomial of ${L}^*(K=I)$ or $L(K=I)$, which is
$f(\lambda)={\lambda}^3-{\lambda}^2+2\lambda-1$.

{\bf Definition 3.} We say that $\lambda$ is an unremovable spectrum
of system (1) with state feedback $u=Kx$,  if there exists $0\ne
X\in {\cal S}_n$, such that for any $K\in {\cal R}^{m\times n}$,
\begin{equation}
X(A+BK)+(A+BK)'X+(C+DK)'X(C+DK)=\lambda X
\label {eq 2.8}
\end{equation}
holds.

Below, we give a theorem with respect to the unremovable spectrum.

{\bf Theorem 2 \footnote {This theorem was first pointed out as a
conjecture  in `` W. Zhang,  Study on algebraic Riccati equation
arising from infinite horizon stochastic LQ optimal control,  Ph.D
Thesis, Hangzhou: Zhejiang University, 1998."}.} $\lambda$ is an
unremovable spectrum of system (1) iff there exists $0\ne X\in {\cal
S}_n$, such that the following three equalities
\begin{equation}
XA+A'X+C'XC=\lambda X, \ \ XB+C'XD=0, \ \ D'XD=0
\label {eq 09}
\end{equation}
hold.

{\bf Proof.} Note that (\ref {eq 2.8}) can be written as
\begin{equation}
XA+A'X+C'XC+(XB+C'XD)K+K'(XB+C'XD)'+K'D'XDK=\lambda X, \label {eq
010}
\end{equation}
so if (\ref {eq 09}) holds, then (\ref {eq 2.8}) automatically
holds, and the sufficiency is proved.

To prove the necessity  of Theorem 2, we first take $K=0$ in (\ref
{eq 2.8}), then
$$
XA+A'X+C'XC=\lambda X
$$
holds. Again, from (\ref {eq 010}), it follows that
\begin{equation}
(XB+C'XD)K+K'(XB+C'XD)'+K'D'XDK=0. \label {eq 011}
\end{equation}
Let $XB+C'XD=F, D'XD=G$, then (\ref {eq 011}) becomes
\begin{equation}
FK+K'F'=-K'GK.
\label {eq 012}
\end{equation}
Since the left hand side of (\ref {eq 012}) is linear with $K$, we
must have $G=0$. In fact, due to the linearity of $K'GK$,
$$
(K+K)'G(K+K)=4K'GK=K'GK+K'GK=2K'GK.
$$
So $K'GK=0$  because of the arbitrarity of $K$, which  is necessary
that $G=0$, i.e., $D'XD=0$. To prove $F=0$ or $XB+C'XD=0$, we note
that (\ref {eq 012}) becomes $K'F'=-FK$.  Denote
$F=(f_{ij})_{n\times m}$, and take
$$
K=K_{ij}=(k_{ls})_{m\times n }=\left\{
\begin{array}{l}
1,\mathrm{\ if\ }l=i, s=j,\\
0\mathrm{,\ otherwise.}
\end{array}
\right.
$$
From $K'F'=-FK$, one knows $f_{1i}=f_{2i}=\cdots=f_{ni}=0$. Set
$i=1,2,\cdots, n$, then $f_{ij}=0$ for $i=1,2,\cdots,n,
j=1,2,\cdots, m$, that is $F=0$.  The proof of Theorem 2 is
completed.

{\bf Remark 2.}  Commonly,  we call  $\lambda$   an unremovable
spectrum of (\ref{eq 1.3}) with state feedback $u=Kx$, if for any
feedback gain $K$, $\lambda \in \sigma(A+BK)$. We should point out
that, when $C=0, D=0$, Definition 3 coincides with  the conventional
one  of system (\ref{eq 1.3})  illustrated by the following
Proposition, which may be new even in linear system  theory.

{\bf Proposition 1.}  For $C=0, D=0$, the following conditions are
equivalent:

1)\ There exists $0\ne X\in {\cal S}_n$ satisfying (\ref{eq 09}).

2)\  There exists a nonzero complex vector
$\xi$ satisfying (\ref{eq quo}).

3)\  For any $K\in {\cal R}^{m\times n}$, $\lambda\in \sigma(A+BK)$.

{\bf Proof.} 1)$\Leftrightarrow$ 2). If 1) holds, then there exists
$0\ne X\in {\cal S}_n$ satisfying
$$
XA+A'X=\lambda X, \ \ XB=0,
$$
which is equivalent to
\begin{equation}
XA+A'X=\lambda X, \ \ XBB'=0.
\label {eq prop}
\end{equation}
(\ref{eq prop}) can be written as
\begin{equation}
(I\otimes A'+A'\otimes I)\vec{X}=\lambda \vec{X}, \ \ (I\otimes
BB')\vec{X}=0, \label {eq kro}
\end{equation}
where $\vec {X}$ denotes the vector formed by stacking the rows of $X$ into one long vector, i.e.,
$$
\vec{X}=[X_{11}, X_{12},\cdots, X_{1n}, X_{21}, X_{22}, \cdots, X_{2n},\cdots,
X_{n1}, X_{n2},\cdots, X_{nn}]',
$$
$F\otimes G$ denotes the Kronecker product of  two matrices $F$ and
$G$. Applying the above two facts, one can immediately shows
1)$\Rightarrow$ 2).

Fact 1 (Ortega, 1987). If ${\xi}_1, {\xi}_2,
\cdots,
{\xi}_p$ are linearly independent eigenvectors of $A'$,
then the Kronecker product
$$
{\xi}_i\otimes {\xi}_j,\ \ i,j=1,2,\cdots, p
$$
are linearly independent  eigenvectors of
$(I\otimes A'+A'\otimes I)$.

Fact 2 (Ortega, 1987). If $e_1, e_2,\cdots, e_n$ are the coordinate vectors of ${\cal R}^n$,
$\zeta _1,
\zeta _2, \cdots, \zeta _q$ are linearly independent eigenvectors of $BB'$ corresponding to
zero eigenvalue , then
the Kronecker product of $e_i$ and $\zeta _j$
$$
e_i\otimes \zeta _j, \ \  i=1,2,\cdots, n, j=1, \cdots, q
$$
are linearly independent eigenvectors of $I\otimes BB'$  corresponding to zero eigenvalue.

2)$\Rightarrow$ 1) is simple, we only need to take $X=\xi{\xi}'$.  So 1)$\Leftrightarrow$ 2).

2)$\Leftrightarrow$ 3). From 2), it is easy to show for any $K$, $\lambda\in \sigma(A'+K'B')
=\sigma(A+BK)$, so 2)$\Rightarrow$ 3).

3)$\Rightarrow$ 2). If $Re (\lambda)\ge 0$, and $\lambda\in
\sigma(A+BK)$ with any $K$, but 2) does not hold, then  by PBH
criterion , $(A,B)$ is stabilizable. So there exists $K_1$ such that
$\sigma(A+BK_1)\subset {\cal C}^-$, accordingly, it is impossible
that $\lambda\in \sigma(A+BK)$ with $Re (\lambda)\ge 0$. If
$Re(\lambda)<0$, then $\mu=-\lambda\in \sigma(-A-BK)$ with
$Re(\mu)>0$. By the same discussion, we can still show that 3)
implies 2). Therefore,  2)$\Leftrightarrow$ 3).  Proposition 1 is
proved.

By Theorem 2 and Proposition 1, deterministic PBH criterion can be
stated in another form as follows:

{\bf Corollary 1.} When $C=0, D=0$ in (1), $(A,B)$ is stabilizable
iff all the  unremovable spectra of system (\ref{eq 1.3}) belong to
${\cal C}^-$, that is, there does not exist a nonzero $X \in {\cal
S}_n$, $Re(\lambda)\ge 0$,  satisfying
\begin{equation}
XA+A'X=\lambda X,\ \  XB=0. \label {eq xavv}
\end{equation}

From Corollary 1, it is natural to conjecture that system (1) is
stabilizable iff  its all unremovable spectra in ${\cal C}^-$, or
there does not exist $0\ne X\in {\cal S}^n$, $Re(\lambda)\ge 0$
satisfying (\ref {eq 09}). Unfortunately, the following example
shows  that  this conjecture is  not true, so PBH criterion  cannot
be generalized to the stabilizability of stochastic systems, which
reveals the essential difference between deterministic and
stochastic systems with control entering into diffusion term.

{\bf Example 2.} Consider  a one-dimensional case of system (1).
Take  $D\ne 0$, then there does not exist  $X\ne 0$, $Re(\lambda)\ge
0$, satisfying (\ref {eq 09}).  But from Theorem 1, system (1) is
stabilizable iff  $B^2+2BCD-2AD^2>0$. Obviously, $B^2+2BCD-2AD^2>0$
is not equivalent to $D\ne 0$.

But can we be sure that stochastic PBH criterion holds with $D=0$?
The answer is still no;  see the following example:

{\bf Example 3.}  In system (1), take $D=0$,
$$
A=\left[
\begin{array}{cc}
0 &  1\\
0 &  -1
\end{array}
\right],  \quad
B=\left[
\begin{array}{cc}
0 \\
1
\end{array}
\right], \quad
C=\left[
\begin{array}{cc}
1 & -1\\
0 &  0
\end{array}
\right],
$$
$$
Q=\left[
\begin{array}{cc}
1 &  0\\
0 &  1
\end{array}
\right],  \quad
R=\left[
\begin{array}{cc}
1 & 0 \\
0  & 1
\end{array}
\right].
$$
One can test that  for any $X\in {\cal S}_2$, the following equation
$$
XA+A'X+C'XC=\lambda X, \ \ XB=0, \ \ Re(\lambda)\ge 0
$$
does not have nonzero solution $X$, but $(A,B;C,0)$ is not
stabilizable. Because by Zhang (2000), $(A,B;C,0)$ is stabilizable
iff GARE
\begin{equation}
PA+A'P+C'PC-PBR^{-1}B'P+Q=0
\label {eq lz}
\end{equation}
has a positive definite solution $P$, but for our given data, the
solutions of GARE (\ref{eq lz}) is
$$
P_1=\left[
\begin{array}{cc}
-1 &  0\\
0 &  0
\end{array}
\right],  \quad
P_2=\left[
\begin{array}{cc}
-1 & 0 \\
0  & -2
\end{array}
\right].
$$
When does the PBH criterion hold for stochastic stabilizability? For
the special case, we obtain the following theorem.

{\bf Theorem   3 (Stochastic PBH criterion).} For  $D=0$, if there
exists a matrix $C_1$, such that for any $X\in {\cal S}_2$,
 $C'XC= XC_1+C'_1X$,   then system (1)  is stabilizable
iff  its all unremovable spectra belong to ${\cal C}^-$, i.e., there
does not exist nonzero  $X\in {\cal S}^n$, such that
\begin{equation}
XA+A'X+C'XC=\lambda X, XB=0. Re(\lambda)\ge 0. \label {eq st}
\end{equation}

{\bf Proof.}
By Theorem 1 and Remark 1, system (1) is stabilizable iff there exists an
$K\in {\cal R}^{m\times n}$,
$\sigma ({\cal L}^*_K)\subset {\cal C}^-$. Note that for $X\in {\cal S}_n$,
\begin{eqnarray}
{\cal L}^*_K (X)&=&X(A+BK)+(A+BK)'X+C'XC \nonumber \\
&=& X(A+C_1+BK)+(A+BK+C_1X)'X. \nonumber
\end{eqnarray}
We know that system (1) is stabilizable iff the following
deterministic linear system
$$
\dot {z}=(A+C_1)z+Bu
$$
is stabilizable, which, from Corollary 1, is equivalent to that
there does not exist nonzero $X$ satisfying
$$
X(A+C_1)+(A+C_1)'X=XA+A'X+C'XC=\lambda X, \ \ XB=0.
$$
Theorem 3 is proved.

Another problem is on the spectral  placement. It is well known that
in deterministic systems, a necessary and sufficient condition for
complete controllability of $(A,B)$ is that the spectrum of (\ref{eq
1.3}) can be arbitrarily assigned ( De Carlo 1989). But for
stochastic systems, it is hard to relate  controllability with
spectrum in many existing definitions of stochastic controllability;
see Peng (1994), Bashirov and Kerimov (1997), Mahmudov (2001). To
illustrate the complexity, consider system (1) with ungiven initial
state, on which the exactly-terminal controllability was introduced
by Peng (1994) as follows:

{\bf Definition 4. } System (1) is called exactly
terminal-controllable, if for any $\xi\in L^2 (\Omega,{\cal F},{\cal
P};{\cal F}_T)$, there exists at least one admissible control
$u(t)\in L^2(\Omega,{\cal F},{\cal P};{\cal F}_t)$, and initial
state $x_0\in {\cal R}^n$, such that the corresponding trajectory
satisfies $ x(T)=\xi$, where $L^2 (\Omega,{\cal F},{\cal P};{\cal
F}_T)$ denotes  all ${\cal F}_t$-adapted, measurable and square
integrable processes.

The following example shows that even if system (1) is
exactly-terminal controllable, its spectrum cannot be assigned
arbitrarily.

{\bf Example 4.} We still consider one-dimensional case.  Assume (1)
is exactly-terminal controllable, then $D\ne 0$ (Peng, 1994). By a
simple computation, for any $K$,
$$
\sigma ({\cal L}^*_K)=\lambda=D^2K^2+2(CD+B)K+2A+C^2\ge 2A-\frac
{(CD+B)^2}{D^2}.
$$
So the spectrum cannot take arbitrary value in ${\cal R}$.

We note that (\ref {eq 05}) can also be written as
\begin{equation}
\dot {\vec{X}}= [(A+BK)\otimes
I+I\otimes(A+BK)+(C+DK)\otimes(C+DK)]\vec{X}:=L^0(K)\vec{X}. \label
{eq kro}
\end{equation}
 Since
$\lim_{t\to\infty} E[x(t)x'(t)]=0 \Leftrightarrow \lim_{t\to\infty}
\vec{X}=0$,  $(A,B;C,D)$ is stabilizable iff there exists a $K$,
such that $\sigma (L^0(K))\subset {\cal C}^-$. Kleinman (1969)
further asserted that  $L^0(K)$ will have several repeated
eigenvalues, this assertion is not true; see Zhang (1998). While
each of $L^0(K)$ and $L(K)$ describes the stabilizability of system
(1), it is natural to ask how relation with the spectrum between
$L^0(K)$ and $L(K)$? In general, we have the following result.

{\bf Proposition 2.}

1) If $\lambda\in \sigma (L(K))$, then $\lambda\in \sigma (L^0(K))$;

2) If $\lambda \in \sigma (L^0(K))$, but $\lambda$ does not belong
to $\sigma (L(K))$, then there exists $X\in {\cal C}^{n\times n}$,
such that $L^0(K)\vec{X}=\lambda \vec{X}$, $X=-X'$.

{\bf Proof.}  If we Define   linear operator ${\cal L}^0_K $ as
follows:
$$
{\cal L}^0_K : X\in {\cal C}^{n\times n}\longmapsto (A+BK)X+X(A+BK)'+(C+DK)X(C+DK)'
$$
then it is easily seen $\sigma (L^0(K))=\sigma({\cal L}^0_K )$. If $\lambda\in \sigma (L(K))$,
then there exists $X\in {\cal S}_n$, such that ${\cal L}_K(X)=\lambda X$, therefore,
${\cal L}^0_K (X)=\lambda X$, i.e, $\lambda \in \sigma (L^0(K))$.

As to prove 2), we note that if $\lambda \in \sigma (L^0(K))$, then there exists
 $X\in {\cal C}^{n\times n}$ satisfying
${\cal L}^0_K (X)=\lambda X$. By symmetry, there also has ${\cal
L}^0_K (X')=\lambda X'$. If $X\ne -X'$, then from
$$
{\cal L}^0_K (X+X')={\cal L}^0_K (X)+{\cal L}^0_K (X')=\lambda X+\lambda X'=\lambda(X+X'),
$$
$(X+X')\ne 0\in {\cal S}_n$ is an eigenvector with respect to the
eigenvalue $\lambda$.  Therefore,  $\lambda\in \sigma (L(K))$, which
results in a contradiction, and,  accordingly,  $X=-X'$.

Theorem 1 or Proposition 1 has theoretical value in studying the
spectral  allocation and the solutions of GARE.  Obviously, from the
practical point of view, both  of  them are  not convenient for
testing the stabilizability of $(A,B;C,D)$.  Ait Rami and Zhou
(2000) presented an efficient method expressed by an LMI.

\section{ The spectral characterization for stochastic observability and detectability}

In this section, we apply the spectral technique  to study the
observability and detectability of system (1). We first give a
spectral criterion for exact observability of stochastic systems,
and  then by duality, we define stochastic detecterbility.  All
these concepts play  critical roles in many fields,   and  some
applicable examples can be found in section 5 of this paper.

{\bf Definition 5:} Consider the following stochastic system with
measurement equation:
\begin{equation}
\left\{
\begin{array}{l}
dx(t)=(Ax(t)+Bu(t))\,dt+(Cx+Du)\,dw(t), x(0)=x_0\in {\cal R}^n,\nonumber\\
y(t)=Qx(t). \label {eq de5}
\end{array}
\right.
\end{equation}
We call  $x_0\in {\cal R}^n $ an unobservable state, if let
$u(t)\equiv 0$, then,  for any $T>0$, the output response with
respect to $x_0$,   is always equal to zero, i.e.,
$$
y(t)\equiv 0, \quad a.s. , \ \  t\in [0,T].
$$

{\bf Definition 6:} System (\ref{eq de5}) is called exactly observable, if there is no unobservable
initial state (except zero initial state).

{\bf Remark 3.} Definition 6 can also be equivalently expressed as
follows: We call system (\ref{eq de5}) exactly observable, if for
arbitrary $x_0\in {\cal R}^n$, $x_0\ne 0$, there exists a  $t>0$
such that $y(t)=Qx(t)\ne 0$, where $x(t)$ is the solution to
stochastic differential equation (SDE)
\begin{equation}
dx=Ax\,dt +Cx\,dw, x(0)=x_0\in {\cal R}^n.
\label {eq obz}
\end{equation}

For simplicity, when system (\ref{eq de5}) is exactly observable, we
also call  $[A, C|Q]$ exactly observable.

{\bf Theorem 4.} $[A, C|Q]$ is exactly observable iff there does not
exist nonzero $X\in {\cal S}_n$, such that
\begin{equation}
XA'+AX+CXC'=\lambda X,\ \  QX=0, \ \ \lambda\in {\cal C}. \label{eq
the}
\end{equation}

{\bf Proof.} Let $X(t)=E[x(t)x'(t)], Y(t)=E[y(t)y'(t)]$, where
$x(t)$ is the solution of (\ref{eq obz}). As in the previous
discussion, $X(t)$ satisfies
\begin{equation}
\dot {X}(t)=AX+XA'+CXC':={\cal L}_{A,C}(X(t)), \  X(0)=x_0x'_0.
\label {eq obz1}
\end{equation}
By Remark 3, we know that  $[A, C|Q]$ is exactly observable iff for
arbitrary $X_0=x_0x'_0\ne 0$, there exists a  $t>0$ such that
\begin{equation}
Y(t)=Ey(t)y'(t)=QX(t)Q'\ne 0.
\label {eq obz2}
\end{equation}
From the proof of Theorem 3, (\ref{eq obz1}) is equivalent to
\begin{equation}
\dot {\tilde {X}} =L(A,C){\tilde X},
\label {eq obz3}
\end{equation}
where $L(A,C)$ is the induced  matrix by  operator ${\cal L}_{A,C}$.
Secondly, due to $X(t)\ge 0$ for all $t\ge 0$,  (\ref{eq obz2}) is
equivalent to
\begin{equation}
Y_1=QX(t)\ne 0, \label{eq obz5}
\end{equation}
which is equivalent to
\begin{equation}
\vec {Y_1}=\vec {L}(Q) \tilde X\ne 0. \label {eq obz6}
\end{equation}
${\vec { L}}(Q)$ is one $n^2\times \frac {n(n+1)} 2$-order matrix, which is uniquely determined
by $Q$. So (\ref{eq de5}) is exactly observable iff
the deterministic system
\begin{equation}
\left\{
\begin{array}{l}
\dot {\tilde {X}} =L(A,C){\tilde X},\nonumber\\
\vec {Y_1}={\vec { L}}(Q) \tilde X
\label {eq de6}
\end{array}
\right.
\end{equation}
is completely observable. By PBH  criterion for observability,
(\ref{eq de6}) is completely observable iff there does not exist
$0\ne \xi\in {\cal C}^{\frac {n(n+1)}2}$, such that
$$
L(A,C)\xi=\lambda \xi, {\vec { L}}(Q)\xi=0, \lambda\in {\cal C}.
$$
By our definition with $L(A,C),{\vec { L}}(Q)$,   which is
equivalent to that there does not exist $0\ne X\in {\cal S}_n$
satisfying (\ref{eq the}), the proof of Theorem 4 is completed.

We give the following example to illustrate the notion of ${\vec { L}}$.

{\bf Example 5.} In (\ref{eq obz5}), taking
$$
Y_1=\left[
\begin{array}{cc}
y^1_{11} &  y^1_{12}\\
y^1_{21} &  y^1_{22}
\end{array}
\right],  \quad
Q=\left[
\begin{array}{cc}
3 & 2\\
5 &  7
\end{array}
\right], \quad
X=\left[
\begin{array}{cc}
x_{11} & x_{12}\\
x_{12} &  x_{22}
\end{array}
\right].
$$
From $Y_1=QX$, we have
\begin{equation}
\left\{
\begin{array}{l}
y^1_{11}=3x_{11}+2x_{12},\nonumber\\
y^1_{12}=3x_{12}+2x_{22},\nonumber\\
y^1_{21}=5x_{11}+7x_{12},\nonumber\\
y^1_{22}=5x_{12}+7x_{22}.
\label {eq exam}
\end{array}
\right.
\end{equation}
(\ref{eq exam}) can be written in the matrix form as
$$
\vec {Y_1}=\left[
\begin{array}{ccc}
3 & 2 & 0\\
0 & 3 & 2\\
5 & 7  & 0\\
0 & 5 & 7
\end{array}
\right]\tilde X:= {\vec { L}}(Q) \tilde X.
$$
So
$$
{\vec { L}}(Q)=\left[
\begin{array}{ccc}
3 & 2 & 0\\
0 & 3 & 2\\
5 & 7  & 0\\
0 & 5 & 7
\end{array}
\right].
$$
{\bf Remark 4.} Liu (1999) proved a dual principle, which asserted
that  $[A,C|Q]$ is exactly observable  iff the following system is
exactly controllable.
$$
-dx=(A'x+C'z+Q'u)\,dt-zdw.
$$
{\bf Corollary 2.} If there does  not exists nonzero $X\in {\cal
S}_n$ satisfying
\begin{equation}
XA'+AX=\lambda X,\ \  QX=0, \ \ \lambda\in {\cal C}, \label {eq co5}
\end{equation}
then for any real matrix $C$ of suitable  dimension,
\begin{equation}
XA'+AX+CXC'=\lambda X,\ \  QX=0, \ \ \lambda\in {\cal C}
\label {eq co3}
\end{equation}
does  not have nonzero solution $X\in {\cal S}_n$.

{\bf Proof.}  By the same discussion as in Theorem 4, we can prove
that  $(Q,A)$ is completely observable iff there does not exist
nonzero $X\in {\cal S}_n$ satisfying (\ref{eq co5}). Secondly, from
Liu (1999), we know that $[A,C|Q]$ is exactly observable iff $ Rank
(P_0)=n $, where
$$ P_0=[ Q', A'Q', C'Q', A'C'Q', C'A'Q', (A')^2Q', (C')^2Q',\cdots]'.
$$
Obviously, $Rank(Q',A'Q', (A')^2Q',\cdots, (A')^{n-1}Q')'\le
Rank(P_0)$, so if $(Q,A)$ is completely observable, then so does
$[A,C|Q]$. Applying Theorem 4, Corollary 2 is immediately derived.

The following proposition may be useful in some cases, it is a
generalized version of Lemma 4.1 of Wonham (1968).

{\bf Proposition 3.} Let
\begin{equation}
Q'Q+D'D=F'F
\label {eq won}
\end{equation}
and $G_1, G_2$ are any real matrices of suitable dimension. If
$[A,C|Q]$ is exactly observable, then so does $[A+G_1D, C+G_2D|F]$.

{\bf Proof.} If $[A+G_1D, C+G_2D|F]$  is not exactly observable,
then there exists $0\ne X\in {\cal S}_n$, such that
$$
X(A+G_1D)'+(A+G_1D)X+(C+G_2D)X(C+G_2D)'=\lambda X, \ \  FX=0, \ \
\lambda\in {\cal C}.
$$
From $FX=0$ and (\ref{eq won}), we have $QX=0, DX=0$, together with the above equations,
one has
$$
XA'+AX+CXC'=\lambda X,\ \ QX=0, \ \ \lambda\in {\cal C}.
$$
By Theorem 4, $[A,C|Q]$ is not exactly observable, which
contradicts the given conditions.

Based on stabilizability, we can define stochastic detectability via duality.

{\bf Definition 7.} We say that  $[A, C|Q]$ is stochastic
detectable, if $(A', Q'; C',0)$ is stabilizable.

From Theorem 3, we have

{\bf Corollary 3.} If  $[A, C|Q]$  is  stochastic detectable,  then
there does not exist nonzero $ X\in {\cal S}_n$, such that
\begin{equation}
XA'+AX+CXC'=\lambda X, QX=0, Re(\lambda)\ge 0. \label {eq co2}
\end{equation}
If there exists $C_1$  such that for any $X\in {\cal S}_n$,
$CXC'=XC'_1+C_1X$, then (\ref{eq co2}) is also a necessary condition
for stochastic detectability.

There does not have any implication  between exact observability and
stochastic detectability; see the following examples:

{\bf Example 6.}  Take
$$
A=\left[
\begin{array}{cc}
0 &  0\\
1 &  -1
\end{array}
\right],  \quad
Q=\left[
\begin{array}{cc}
0 & 1
\end{array}
\right], \quad
C=\left[
\begin{array}{cc}
1 & 0\\
-1 &  0
\end{array}
\right].
$$
One can easily test that (\ref{eq the}) does not have nonzero
solution $X\in {\cal S}_2$, so $[A,C|Q]$ is exactly observable.
However,  $[A,C|Q]$ is not stochastic detectable, because $(A',Q',
C',0)$ is not stabilizable as shown in Example 3.

{\bf Example 7.} Take $Q=0$,
$$
A=\left[
\begin{array}{cc}
-1 &  0\\
0 &  -1
\end{array}
\right], \quad
C=\left[
\begin{array}{cc}
-1 &  0\\
0 &  -1
\end{array}
\right],
$$
then it is easily tested that $[A,C|Q]$ is stochastic detectable,
but it is not exactly observable.

{\bf Proposition 4. } If  $[A,C|Q]$ is stochastic  detectable, then
$(Q,A)$ is detectable.

{\bf Proof.} By Definition 7, if $[A,C|Q]$ is  exactly detectable, then  there exists a
constant matrix
$H'$, such that
\begin{equation}
dx=(A'+Q'H')x\,dt +C'x\,dw
\end{equation}
is asymptotically  mean square stable, which implies (Has'minskii,
1980) that
\begin{equation}
dx=(A'+Q'H')x\,dt
\end{equation}
is asymptotically stable, so $(Q,A)$ is detectable.

\section { On weak stabilizability of stochastic systems}

In this section, we study the weak stabilizability of system (1),
which describes a class of weak  stability, and has close
relationship with strong solutions of GAREs.

{\bf Definition 8.} We say system (1) is weakly stabilizable, if
there exists a matrix $K\in {\cal R}^{m\times n}$, via the state
feedback $u(t)=Kx(t)$, the closed-loop system (\ref{eq 78}) is
two-stable (Has'minskii, 1980), i.e., for each $\varepsilon>0$,
there exist an $\delta>0$, such that
$$
E|x(t,x_0)|^2<\varepsilon
$$
whenever $t\ge 0$ and $|x_0|<\delta$.

Our result  is given as follows:

{\bf Proposition  5.}
 System (1) is weakly stabilizable iff  one of the following conditions holds:

1)\  There exists $K\in {\cal R}^{m\times n}$, such that
$\sigma (L(K))=\{\lambda _i, i=1,2,\cdots, n(n+1)/2\}\subset {\cal C}^{-,0}$,
and whenever $Re\lambda _i=0$, all associated Jordan blocks of $\lambda _i=0$
are $1\times 1$;

2)\ There exists $P>0$, such that
\begin{equation}
PL(K)+L'(K)P\le 0  .
\label {eq pl}
\end{equation}
{\bf Proof.} By Definition 8, system (1) is weakly stabilizable iff
there exists $K$ with that system (\ref {eq 78}) is   weakly stable.
As done in Theorem 1, noting that for each  $\varepsilon>0$, there
exists a $\delta>0$, such that
$$
E|x(t,x_0)|^2<\varepsilon
$$
whenever $t\ge 0$ and $|x_0|<\delta$ , which  is equivalent to that
for each  $\varepsilon>0$, there exists a $\delta>0$, such that
$$
||\tilde {X}(t,x_0)||<\varepsilon
$$
whenever $t\ge 0$ and $||\tilde {X}_0||<\delta$. From (\ref {eq 06}) and
Theorem 5.2.3 of Ortega (1987), the
latter is equivalent to that  $L(K)$ is weakly negative stable, which completes the proof of 1).

2) is a simple corollary of Theorem 5.4.3 of Ortega (1987).

The following theorem is a sufficient condition for weak
stabilizability expressed  by LMIs and Lyapunov-type inequalities,
which will be used later. Analogous results for stabilizability can
be found in Ait Rami and Zhou (2000).

{\bf Theorem 5.} System (1) is weakly stabilizable if one of the following conditions holds.

1) There exist $K\in {\cal R}^{m\times n}$ and $P>0$, such that
\begin{equation}
P(A+BK)+(A+BK)'P+(C+DK)'P(C+DK)\le 0.
\end{equation}

2) There exist $K\in {\cal R}^{m\times n}$ and $P>0$, such that
\begin{equation}
(A+BK)P+P(A+BK)'+(C+DK)P(C+DK)'\le 0.
\label {eq 2)}
\end{equation}

3) There exist  matrices $P$ and $Y$, such that the following LMI
 holds.
\begin{equation}
\left[
\begin{array}{cc}
 AP+PA'+BY+Y'B' &  CP+DY \\
  PC'+Y'D'  &   -P
\end{array}
\right]\le 0 .
\label {eq cd}
\end{equation}

{\bf Proof.}
 If 1) holds, then by Dynkin's formula (Oksendal,1998), we have
\begin{eqnarray}
&{}&Ex'(t)Px(t)=x'_0Px_0+E\int^t _0 {\cal A} (x'(s)Px(s))\,ds=x'_0Px_0 \nonumber \\
&+&  E\int^t _0 (<(A+BK)x,\frac {\partial(x'PX)}{\partial x}>+
\frac 1 2<(C+DK)x, \frac {{\partial}^2(x'PX)}{\partial
x^2}(C+DK)x>)\,ds \nonumber \\
&=&x'_0Px_0
+ E\int^t _0 x'(P(A+BK)+(A+BK)'P+(C+DK)'P(C+DK))x\,ds\nonumber \\
&\le & x'_0Px_0, \label{eq 18}
\end{eqnarray}
where ${\cal A}$ is the infinitesimal generator of $x(t)$, i.e.,
the trajectory of stochastic system (\ref {eq 78}). From (\ref{eq
18}), for each $\varepsilon>0$,
\begin{eqnarray}
E|x(t)|^2 &\le& \max (1, \lambda _{max} (P))E|x(t)|^2\nonumber \\
&\le& \max (1, \lambda _{max} (P))
\frac {\lambda _{max}(P)}{\lambda _{min} (P)} |x_0|^2\nonumber \\
&=& \max (\frac {\lambda _{max}(P)}{\lambda _{min} (P)},
\frac {{\lambda}^2 _{max}(P)}{\lambda _{min} (P)}) |x_0|^2 :=C_0 |x_0|^2<\varepsilon
\end{eqnarray}
whenever $|x_0|<\delta:=\frac {{\varepsilon}^{1/2}}{{C_0}^{1/2}}$.
So system (1) is weakly stabilizable.

If 2) holds,  by the same discussion as in 1), we can prove that
the dual system of (\ref {eq 78})
\begin{equation}
dx(t)=(A+BK)'x(t)\,dt+(C+DK)'x(t)\,dw(t), x(0)=x_0
\end{equation}
is weakly stable. From Theorem 1 and Remark 1, this is equivalent to the weak stabilizability of
(\ref {eq 78}).

As to that 3) implies (1) stabilizable, this is due to the
equivalence of 2) and  3). Set $Y=KP, P>0$, then (\ref{eq 2)}) can
be written as
\begin {equation}
AP+BY+PA'+Y'B'+(CP+DY)P^{-1}(CP+DY)'\le 0. \label {eq nn}
\end {equation}
Applying Schur's lemma, (\ref {eq nn}) is equivalent to (\ref{eq
cd}),  which ends our proof.

{\bf Remark 5.} From Proposition 5, we have every reason to
conjecture that Theorem 5 should be not only a sufficient, but also
a necessary condition for weak stabilizability, but how do we prove
its true?

{\bf Remark 6.} Taking $P=I$ in (\ref{eq 18}), apparently,
$(A,B;C,D)$ being weak stabilizable implies that $(A,B)$ is weakly
stabilizable.

The following Proposition will be used in Section 5.

{\bf Proposition 6.} For any real matrix $K\in {\cal R}^{m\times n}$,
if $\sigma ({\cal L}^*_K)\subset {\cal C}^{-,0}$, then $\sigma (A+BK)\subset {\cal C}^{-,0}$.

{\bf Proof.}  Since  $\sigma ({\cal L}^*_K)\subset {\cal C}^{-,0}$, so for any $\varepsilon>0$,
$\sigma ({\cal L}^{*,\varepsilon }_K)\subset {\cal C}^-$, where
${\cal L}^{*,\varepsilon }_K (\cdot)$ is defined by
$$
{\cal L}^{*,\varepsilon} _K: X\in {\cal S}^n\longmapsto X(A-\varepsilon I+BK)+
(A-\varepsilon I+BK)'X+(C+DK)X(C+DK)'.
$$
By Theorem 1 and Remark 1, the stochastic system
$$
dx(t)=(A-\varepsilon I+BK)x\,dt+(C+DK)x\,dw(t), x(0)=x_0\in {\cal
R}^n
$$
is stable in mean square sense, which implies (Has'minskii, 1980)
$\sigma (A-\varepsilon I+BK)\subset {\cal C}^-$. Let $\varepsilon\to
0$, by the continuity of spectrum ( Sontag,1990), we have $\sigma
(A+BK)\subset {\cal C}^{-,0}$.

We should point out that the inverse of Proposition 6 does not hold
even if $\sigma (A+BK)\subset {\cal C}^-$;  see the following
example.

{\bf Example 8.} Take  $A=-I, C=\left[
\begin{array}{cc}
3 & 0\\
0 &  4
\end{array}
\right], K=0$, $B$  and  $D$ are arbitrary, then $\sigma
(A+BK)=\sigma(A)=\{-1,-1\}$, but one can compute $\sigma ({\cal
L}^*_{K=0})=\{7,10,14\}$.

\section{ Some applications}

\subsection {Applications to GAREs}

In this section we apply spectrum technique to study GARE
\begin{equation}
\left\{
\begin{array}{l}
PA+A'P+C'PC+Q-(PB+C'PD)(R+D'PD)^{-1}(B'P+D'PC)=0,\nonumber\\
R+D'PD>0.
\label{eq 4.1}
\end{array}
\right.
\end{equation}
This equation has many applications in infinite horizon linear
quadratic optimal control, stochastic stability, filtering, etc.,
see Bensoussan (1982,1992), Liu (1999), Gao and Ahmed (1987), Wonham
(1968), Zhang (1998),  and the references therein. Equation (\ref{eq
4.1}) is a generalized version of deterministic algebraic Riccati
equation (DARE)
\begin{equation}
PA+A'P+Q-PBR^{-1}B'P=0.
\label {eq 4.2}
\end{equation}
It is well known that if $P\in {\cal S}^n$ is a solution of (\ref{eq
4.2}), and $\sigma (A-BR^{-1}B'P)\subset {\cal C}^-$, then $P$ is
called a feedback stabilizing solution; If $\sigma
(A-BR^{-1}B'P)\subset {\cal C}^{-,0}$, then $P$ is called a strong
solution (Park and Kailath 1997). In other words, the classification
of solutions of algebraic Riccati equations is according to the
stability of the closed-loop system (taking $u=-R^{-1}B'Px$)
$$
\dot x=(A-BR^{-1}B'P)x.
$$
But how do we  make a classification for the solutions  of GARE
(\ref{eq 4.1})? Especially, how do we  define its strong solution?
We note that some authors ( De Souza and Fragoso, 1990) took $P$ as
a strong solution of GARE  (\ref{eq 4.1}), if the spectral set  of
the following deterministic system
\begin{equation}
\dot x=(A-B(R+D'PD)^{-1}(B'P+D'PC))x
\label {eq 4.3}
\end{equation}
satisfying $\sigma (A-B(R+D'PD)^{-1}(B'P+D'PC))\subset {\cal
C}^{-,0}$. However, this definition about strong solution of GARE
(\ref{eq 4.1}) is unreasonable, since from $\sigma
(A-B(R+D'PD)^{-1}(B'P+D'PC))\subset {\cal C}^{-,0}$, we do not know
anything  about the stability of  the closed-loop  system (set the
feedback gain $K=-(R+D'PD)^{-1}(B'P+D'PC)$)
\begin{equation}
dx=(A+BK)x\,dt+(C+DK)x\,dw. \label {eq 4.4}
\end{equation}

Now we give a more reasonable definition  based on spectrum as follows:

{\bf Definition 9.} A solution $P\in {\cal S}_n$ of GARE (\ref {eq
4.1}) is called a feedback stabilizing solution, if $\sigma({\cal
L}^*_{K})\subset {\cal C}^-$ ; $P$ is called a strong solution, if
$\sigma({\cal L}^*_{K})\subset {\cal C}^{-,0}$ , where
$$
K=-(R+D'PD)^{-1}(B'P+D'PC) .
$$

The following theorem is a comparison theorem for GARE (\ref{eq
4.1}), which improves the corresponding results of De Souza and
Fragoso (1990), Ait Rami and Zhou (2000). Firstly, we give a lemma
to be used later.

{\bf Lemma 1 (Ichikawa (1979); Prato and Zabczyk (1992).} The
following Ito-type differential equation
\begin{equation}
dx=Fx+Gx\,dw, \ \ x(0)=x_0
\label {eq 4.00}
\end{equation}
is asymptotically mean square stable, iff for any $Q>0$, the
Lyapunov-type equation
$$
PF+F'P+G'PG=-Q
$$
has a positive solution $P>0$.

When (\ref{eq 4.00}) is asymptotically mean square stable, we also
call  $(F,G)$  stable for short.

{\bf Theorem 6.} Suppose $(A,B;C,D)$ is stabilizable with  the
weighting real matrices $(Q, R)\in {\cal S}_n\times {\cal S}_m$. Let
$\hat{P}$ be any real symmtric solution of the GARE
\begin{equation}
\left\{
\begin{array}{l}
PA+A'P+C'PC+\hat{Q}-(PB+C'PD)(\hat{R}+D'PD)^{-1}(B'P+D'PC)=0,\nonumber\\
\hat {R}+D'PD> 0.
\label {eq 4.5}
\end{array}
\right.
\end{equation}
If $R\ge \hat{R}, Q\ge \hat{Q}$, then GARE (\ref{eq 4.1}) has a
maximal solution $\bar {P}$, $\bar {P}\ge \hat {P}$; Moreover,
$\bar{P}$ is a  strong solution.

{\bf Proof.} Define an operator ${\cal R}$ as follows:
\begin{eqnarray}
&{}&{\cal R}(P,M, N,\varepsilon)\nonumber \\
&=&PA+A'P+C'PC+M +\varepsilon I -(PB+C'PD)(N+D'PD)^{-1}(B'P+D'PC).
\nonumber
\end{eqnarray}
From the given conditions, we know that ${\cal R}(\hat {P},\hat
{Q},\hat {R}, 0)= 0, {\cal R}(\hat {P},Q,R, \varepsilon)>0$,
$R+D'\hat {P}D>0$, then by Theorem 10 and Theorem 12 of Ait Rami and
Zhou (2000), there exist maximal solutions ${\hat {P}}_{max},{\hat
{P}}^0_{max}, {\hat {P}}^\varepsilon_{max}$, respectively, to GARE
(\ref{eq 4.5}), (\ref{eq 4.1}) and GARE
\begin{equation}
\left\{
\begin{array}{l}
PA+A'P+C'PC+Q+\varepsilon I-(PB+C'PD)(R+D'PD)^{-1}(B'P+D'PC)=0,\nonumber\\
R+D'PD>0.
\label {eq 4.6}
\end{array}
\right.
\end{equation}
Moreover, under the constraint of (1), we have
\begin{equation}
V^*(\hat {R},\hat {Q})=
\inf _{u\in {\cal U}_{ad}^\infty}
\{ E\int_0^\infty (x'\hat {Q}x+u'\hat {R}u)\,dt, \lim_{t\to \infty} Ex(t)x'(t)=0\}
=x'_0{\hat {P}}_{max}x_0,
\label {eq 4.7}
\end{equation}
\begin{equation}
V^*(R,Q)=
\inf _{u\in {\cal U}_{ad}^\infty }
\{E\int_0^\infty (x'Qx+u'Ru)\,dt,\lim_{t\to \infty} Ex(t)x'(t)=0\} =
x'_0{\hat {P}}^0_{max}x_0,
\label {eq 4.8}
\end{equation}
\begin{equation}
V^*(R,Q+\varepsilon I)=
\inf _{u\in {\cal U}_{ad}^\infty }
\{E\int_0^\infty (x'(Q+\varepsilon I)x+u'Ru)\,dt,\lim_{t\to \infty} Ex(t)x'(t)=0\}=
x'_0{\hat {P}}^\varepsilon _{max}x_0,
\label {eq 4.9}
\end{equation}
where   ${\cal U}_{ad}^\infty $ denotes all ${\cal F}_t$-adapted,
measurable processes  $u(\cdot): [0,\infty)\times \Omega\longmapsto
{\cal R}^m$ , satisfying
$$
E\int_0^\infty |u(t)|\,dt <\infty.
$$
From (\ref{eq 4.7}), (\ref{eq 4.8}), (\ref{eq 4.9}),  one has
\begin{equation}
{\hat {P}}^\varepsilon _{max}\ge \bar {P}:={\hat {P}}^0_{max}\ge
{\hat {P}}_{max}\ge \hat {P}.  \label {eq 4.0}
\end{equation}
The first part of Theorem 6 is proved.

In what follows, we prove that  $\bar {P}$ is a strong solution.
From (\ref{eq 4.9}), ${\hat {P}}^\varepsilon _{max}$ is monotonic
and  bounded  from below with respect to $\varepsilon$, it is easily
derived
\begin{equation}
\lim_{\varepsilon\to 0} {\hat {P}}^\varepsilon _{max}=\bar {P}.
\label {eq 4.10}
\end{equation}
Now we prove that  ${\hat {P}}^\varepsilon _{max}$ is a feedback
stabilizing solution of GARE (\ref{eq 4.6}). Denote
$A_{\varepsilon}=A+BK_{\varepsilon},
C_{\varepsilon}=C+DK_{\varepsilon}$, ${\hat A}=A+B\hat {K}, \hat
{C}=C+D\hat {K}$, where
$$
K_{\varepsilon}=-(R+D'{\hat {P}}^\varepsilon _{max}D)^{-1}(B'{\hat {P}}^\varepsilon _{max}
+D'{\hat {P}}^\varepsilon _{max}C),
$$
$$
\hat {K}=(R+D'{\hat {P}}D)^{-1}(B'{\hat {P}}+D'{\hat {P}}C).
$$
Note that (\ref {eq 4.5}) and (\ref{eq 4.6}) can  be  respectively
written as
\begin{equation}
\left\{
\begin{array}{l}
{\hat {P}}^\varepsilon _{max}A_{\varepsilon}+A'_{\varepsilon}{\hat {P}}^\varepsilon _{max}
+C'_{\varepsilon}{\hat {P}}^\varepsilon _{max}C_{\varepsilon}+K'_{\varepsilon}RK_{\varepsilon}
+Q+\varepsilon I=0,\nonumber\\
R+D'{\hat {P}}^\varepsilon _{max}D>0, \label {eq 4.11}
\end{array}
\right.
\end{equation}
and
\begin{equation}
\left\{
\begin{array}{l}
{\hat {P}}\hat {A}+{\hat {A}}'{\hat {P}}
+{\hat {C}}'{\hat {P}}\hat {C}+{\hat {K}}'\hat {R} \hat {K}+\hat {Q}=0,\nonumber\\
R+D'{\hat {P}}D>0. \label {eq 4.12}
\end{array}
\right.
\end{equation}
Subtracting  (\ref{eq 4.12}) from (\ref{eq 4.11}), by a series of
computations, we have
\begin{eqnarray}
({\hat {P}}^\varepsilon _{max}-\hat {P})A_{\varepsilon}&+&
A'_{\varepsilon}({\hat {P}}^\varepsilon _{max}-\hat {P})+
C'_{\varepsilon}({\hat {P}}^\varepsilon _{max}-\hat {P})C_{\varepsilon}+\varepsilon I\nonumber\\
&+&(\hat {P}B+C'\hat {P}D)K_{\varepsilon}+(\hat {P}B+C'\hat {P}D)
(\hat {R}+D'\hat {P}D)^{-1}(B'\hat {P}+D'\hat {P}C)\nonumber\\
&+&K'_{\varepsilon}(B'\hat {P}+D'\hat {P}C)+K'_{\varepsilon}(R+D'\hat {P}D)K_{\varepsilon}=0
\label {eq 4.13}
\end{eqnarray}
From (\ref{eq 4.13}), especially noting the last term of the right hand side of the above equation
with
the given condition $R\ge \hat {R}$, we derive
\begin{eqnarray}
({\hat {P}}^\varepsilon _{max}-\hat {P})A_{\varepsilon}&+&
A'_{\varepsilon}({\hat {P}}^\varepsilon _{max}-\hat {P})+
C'_{\varepsilon}({\hat {P}}^\varepsilon _{max}-\hat {P})C_{\varepsilon}\nonumber\\
&\le &-\varepsilon I-(\hat {P}B+C'\hat {P}D)K_{\varepsilon}-(\hat {P}B+C'\hat {P}D)
(\hat {R}+D'\hat {P}D)^{-1}(B'\hat {P}+D'\hat {P}C)\nonumber\\
&-&K'_{\varepsilon}(B'\hat {P}+D'\hat {P}C)-
K'_{\varepsilon}(\hat {R}+D'\hat {P}D)K_{\varepsilon}\nonumber\\
&=& -\varepsilon I-[(\hat {P}B+C'\hat {P}D)+K'_{\varepsilon}(\hat {R}+D'\hat {P}D)]
(\hat {R}+D'\hat {P}D)^{-1}\nonumber\\
&\times &[(\hat {P}B+C'\hat {P}D)+K'_{\varepsilon}(\hat {R}+D'\hat {P}D)]'<0.
\label {eq 4.14}
\end{eqnarray}
We first assert that for any $\varepsilon>0$, ${\hat {P}}^\varepsilon _{max}>\hat {P}$. Otherwise,
by (\ref{eq 4.0}),${\hat {P}}^\varepsilon _{max}-\hat {P}\ge 0$,
$Ker({\hat {P}}^\varepsilon _{max}-\hat {P})\ne \phi$. Let
$0\ne \xi\in Ker({\hat {P}}^\varepsilon _{max}-\hat {P})$, premultiplying by ${\xi}'$ and
postmultiplying by $\xi$ in (\ref{eq 4.14}), then
$$
0\le {\xi}'C'_{\varepsilon}({\hat {P}}^\varepsilon _{max}-\hat
{P})C_{\varepsilon}\xi \le -\varepsilon {\xi}'\xi<0,
$$
which is a contradiction, so ${\hat {P}}^\varepsilon _{max}>\hat
{P}$. Together with (\ref{eq 4.14}),  ${\hat {P}}^\varepsilon
_{max}-\hat {P}$ is a positive solution to the Lyapunov-type
inequality
$$
PA_{\varepsilon}+
A'_{\varepsilon}P+
C'_{\varepsilon}PC_{\varepsilon}<0.
$$
So by Lemma 1, we have $(A_{\varepsilon}, C_{\varepsilon})$ is
stable, i.e., ${\hat {P}}^\varepsilon _{max}$ is a feedback
stabilizing solution. So $\sigma (L^*(K_{\varepsilon}))=\sigma
({\cal L}^*_{K_{\varepsilon}})\subset {\cal C}^-$, while from
(\ref{eq 4.10}), it follows  that
$$
\lim_{\varepsilon\to 0} K_{\varepsilon}
= \bar {K}:=-(R+D'{\bar {P}}D)^{-1}(B'{\bar {P}}+D'{\bar {P}}C).
$$
According to the continuity of spectrum (Lemma A.4.1 of Sontag,
1990), we have $\sigma (L^*(\bar {K}))=\sigma ({\cal L}^*_{\bar
{K}})\subset {\cal C}^{-,0}$. The proof of Theorem 6 is completed.

{\bf Corollary 4.} If $ Q\ge 0, R>0$, system (1) is stabilizable,
then GARE (\ref{eq 4.1}) has a maximal solution, which is also a
strong solution.

{\bf Proof.} Take $\hat {Q}=0, \hat {R} =R$ in (\ref{eq 4.5}), then
(\ref{eq 4.5}) also has a solution $\hat {P}=0$.  Corollary 4 is
immediately deduced from Theorem 6.

{\bf Remark 7.} Under the same condition of Corollary 4, Ait Rami
and Zhou (2000) proved that  GARE (\ref{eq 4.1}) has a maximal
solution. Here, we further assert that  it is also a strong
solution.

{\bf Remark 8.} Taking $D=0$, under the stabilizable condition of $(A,B)$, and
\begin{equation}
{\inf}_K \|\int_0^\infty \exp(t(A+BK)')(C'C)\exp(t(A+BK))\,dt \|<1 ,
\label {eq 4.r}
\end{equation}
De Souza and Fragoso (1990) proved an analogous comparison theorem
to  Theorem 5, which asserted that  (\ref{eq 4.1}) has a solution
$P\ge \hat {P}$ with $\sigma (A+BK)\subset {\cal C}^{-,0}$, where
$K=-R^{-1}B'\bar {P}$. Theorem 6 improves their result because our
given condition is weaker (see Zhang (1998)), but the consequence is
stronger than that of De Souza and Fragoso (1990) from Proposition
6.

{\bf Corollary 5.}  For  $ Q\ge 0, R\ge 0$, if GARE (\ref{eq 4.1})
has a solution $P>0$, then $P$ is a weakly feedback  stabilizing
solution.

{\bf Proof.} Note that GARE (\ref{eq 4.1}) can be written as
\begin{equation}
\left\{
\begin{array}{l}
P(A+BK)+(A+BK)'P
+(C+DK)' P(C+DK)=
-Q-K'RK, \nonumber\\
R+D'PD> 0, K=-(R+D'PD)^{-1}(B'P+D'PC). \label {eq 4.11}
\end{array}
\right.
\end{equation}
By applying 1) of Theorem 5, this corollary is immediately derived.

When does GARE (\ref{eq 4.1}) have a feedback stabilizing  solution?
This problem was discussed in detail in Zhang (1998).

{\bf Example 9.}  Assume all the given data are scalar in (\ref{eq
4.1})  with  $a,b,c,d,q$ and
 $r$
substituting  for $ A,B,C,D,Q$ and $R$, respectively. Taking $q=0,
r=1, d\ne 0$, $(a,b;c,d)$ is stabilizable, then GARE (\ref{eq 4.1})
has two solutions $p_1=0$, $p_2= \frac {2a+c^2}{b^2+2bcd-2ad^2}$.
$p_2$ exists due to the following fact:
 Since $(a,b;c,d)$ is stabilizable,
$d\ne 0$, so  $b^2+2bcd-2ad^2>0$  from Theorem 1. For $p_1=0$, the feedback gain $K_1=0$,
and ${\cal L}^*_{K_1}=2a+c^2$.  By Remark 1, $p_1=0$ is a feedback stabilizing solution, if
$$
\left\{
\begin{array}{l}
b^2+2bcd-2ad^2>0,\nonumber\\
2a+c^2<0.
\end{array}
\right.
$$
By Proposition 5 or direct computation, $p_1=0$ is a  weakly feedback stabilizing solution, if
$$
\left\{
\begin{array}{l}
b^2+2bcd-2ad^2>0,\nonumber\\
2a+c^2=0,
\end{array}
\right.
$$
and $p_1=0$ is a asymptotically feedback   stabilizing solution
(i.e., the closed-loop system is asymptotically stable in
probability; see Has'miskii,1980), if
$$
\left\{
\begin{array}{l}
b^2+2bcd-2ad^2>0,\nonumber\\
2a-c^2<0.
\end{array}
\right.
$$
Analogously, for $p_2$, the corresponding feedback gain $K_2=-\frac{2a+c^2}{b^2+2bcd-2ad^2}$.
$p_2$ is a feedback stabilizing solution, if
$$
\left\{
\begin{array}{l}
b^2+2bcd-2ad^2>0,\nonumber\\
2(a+bK_2)+(c+dK_2)^2<0, \nonumber
\end{array}
\right.
$$
and $p_2$ is a feedback weakly stabilizing solution, if
$$
\left\{
\begin{array}{l}
b^2+2bcd-2ad^2>0,\nonumber\\
2(a+bK_2)+(c+dK_2)^2=0.
\end{array}
\right.
$$
Finally, $p_2$ is a  feedback asymptotically stabilizing solution, if
$$
\left\{
\begin{array}{l}
b^2+2bcd-2ad^2>0,\nonumber\\
2(a+bK_2)-(c+dK_2)^2<0.
\end{array}
\right.
$$

{\bf Remark 9.} From the above example, we can see even if ${\cal
L}^*_{K_i}>0$, the closed-loop system can still achieve some
stochastic stability, which reals the complexity about the study of
solutions of GARE (\ref{eq 4.1}).

\subsection{ Applications to the Lyapunov-type equations}

In this section, we study Lyapunov-type equation
\begin{equation}
PA+A'P+\sum_{i=1}^n C'_iPC_i=-Q.
\label {eq lya1}
\end{equation}
It is well known that a matrix $A$ is Hurwitz or deterministic system
$$
\dot x=Ax, x(0)=x_0
$$
is asymptotically stable
iff the following Lyapunov equation
$$
PA+A'P=-Q, \  Q>0
$$
has a positive definite  solution $P>0$. More generally, if $Q\ge
0$, $(Q^{1/2},A)$ is completely observable, then the above assertion
still holds;  see Jacobson, Martin, Pachter \& Geveci, (1980).  In
the same way, (\ref{eq lya1}) has a close relation with the
stability of the following system
\begin{equation}
dx=Ax\,dt +\sum_{i=1}^n C'_ix\,dw_i, \ \ x(0)=x_0\in {\cal R}^n,
\label{eq equ2}
\end{equation}
where $C_i ,i=1, 2,\cdots, n, $ are $n\times n$ real matrices, $w_i,
i=1,2,\cdots, n,$ are independent standard Wiener processes. In
order to discuss (\ref{eq lya1}), Definitions 6 and 7 should be
modified as follows:

{\bf Definition 10.} We call  $[A,C_1,C_2,\cdots, C_n|Q]$ exactly
observable, if for any nonzero $x_0$, there exists $t>0$ such that
$y(t)=Qx(t)\ne 0$, where $x(t)$ is the solution of (\ref{eq equ2}).

{\bf Definition 11.} We call  $[A,C_1,C_2,\cdots, C_n|Q]$ stochastic
detectable, if there exists a real matrix $H$, such that
$$
dx=(A+HQ)x\,dt +\sum_{i=1}^n C'_ix\,dw_i
$$
is asymptotically  mean square stable.

When (\ref{eq equ2}) is asymptotically mean square stable, we also
call $(A,C_1,C_2,\cdots, C_n)$ stable. In the following context, we
will generalize Lemma 1 to $Q\ge 0$ by means of exact observability
and stochastic detectability.

{\bf Theorem 7.} Assume
$Q\ge 0$, then we have

1) If $[A,C_1,C_2,\cdots,C_n|Q^{1/2}]$ is stochastic  detectable, and Lyapunov-type equation
(\ref{eq lya1})
has a solution $P\ge 0$, then $(A,C_1, C_2,\cdots, C_n)$ is stable.

2) If $[A,C_1,C_2,\cdots,C_n|Q^{1/2}]$ is exactly observable, then $(A,C_1,C_2,\cdots, C_n)$
is stable
iff Lyapunov-type
equation (\ref{eq lya1}) has a positive solution $P>0$.

{\bf The proof of 1).} Consider stochastic system
(\ref{eq equ2}),
by Dynkin's formula, we have
\begin{eqnarray}
0\le Ex'(t)Px(t)&=&x'_0Px_0+E\int_0^t x'(s)(PA+A'P+\sum_{i=1}^n C'_iPC_i)x(s)\,ds \nonumber \\
&=&x'_0Px_0-E\int_0^t x'(s)Qx(s)\,ds, \label {eq dyn}
\end{eqnarray}
so
$$
E\int_0^\infty x'Qx\,ds<\infty.
$$
In addition, since $[A,C|Q^{1/2}]$ is stochastically  detectable,
there exists $H\in {\cal R}^{n\times n}$ such that
\begin{equation}
dx=(A+HQ^{1/2})x(t)dt + Cx(t)\,dw
\label {eq fun}
\end{equation}
is asymptotically  mean square stable, which is also mean square
exponentially stable because of time-invariance (Has'minskii, 1980).
Let $x_H$ be fundamental matrix solution associated with (\ref{eq
fun}), then there exist $\alpha,\beta>0$, such that
$$
E|x_H(t,s)|^2\le \beta exp(-\alpha(t-s), t\ge s.
$$
The solution of (\ref{eq equ2}) can be written as
\begin{equation}
x(t)=x_H(t,0)-\int_0^t x_H(t,s)HQ^{1/2}x(s)\,ds. \label {eq xh}
\end{equation}
By a simple discussion, (\ref{eq xh}) results in $E\int_0^\infty
|x(t)|^2\,dt <\infty$, which derives $\lim_{t\to\infty}
E[x(t)x'(t)]=0$ from Has'minskii (1980), 1) is proved.

The proof of  2): If $(A,C_1,C_2,\cdots, C_n)$ is stable, then
(\ref{eq lya1}) has a solution $P\ge 0$;  see El Ghaoui and Ait Rami
(1996). Now we show $P>0$, otherwise, there exists $x_0\ne 0$, such
that $Px_0=0$. From  (\ref{eq dyn}), for any $T>0$, we have
$$
0\le E\int_0^T x'(s)Qx(s)\,ds = -E[x'(T)Px(T)]\le 0,
$$
which follows $y(t)= Qx(t)\equiv 0, \forall t\in [0,T]$, but this is
impossible because of exact observability, so $P>0$.

If (\ref{eq lya1}) has a positive solution $P>0$, from (\ref{eq
dyn}), we know  that  $V(x(t)):=E[x'(t)Px(t)]$ is monotonically
decreasing and bounded from below with respect to $t$, so
$\lim_{t\to\infty} V(x(t))$ exists. If we let $t_n=nT$, then
\begin{equation}
V(x(t_{n+1}))\le V(x(t))\le V(x(t_{n})),\ \  t_n\le t\le t_{n+1}.
\label {eq mm}
\end{equation}
Again, by (\ref{eq dyn}),
\begin{equation}
V(x(t_{n+1}))- V(x(t_n)) =-E\int_{t_n}^{t_{(n+1)}}
x'(t)Qx(t)\,dt=-E[x'(t_n)Hx(t_n)],
\end{equation}
where $H$ is some positive matrix (Zhang (1998)). Taking limit in
the above, we have
$$
\lim_{n\to\infty} V(x(t_n))=\lim_{n\to\infty} E|x(t_n)|^2=0.
$$
By (\ref {eq mm}), $\lim_{n\to\infty} E[x(t)x'(t)]=\lim_{n\to\infty}
V(x(t))=0,$ so $(A,C_1,C_2,\cdots, C_n)$ is stable.

\section{ Robust stabilization of stochastic systems}

In this section, we study the robust quadratic stabilization of the
following system
\begin{equation}
dx=((A+\Delta A) x+Bu)dt+(Cx+Du)dw,\ \ x(0)=x_0, \label {eq rob1}
\end{equation}
where $\Delta A$ is an uncertain real matrix satisfying the ``matching condition"
$$
\Delta A=EFG, F\in {\cal F}=\{F: F'F\le I, F\in {\cal R}^{k\times
j}\}.
$$
The analogous problems for deterministic systems were  discussed by
Khargonekar, Petersen and Zhou (1990), Petersen (1988), et al. Such
problems have significant sense, since even for stochastic models,
the coefficient matrices is not necessary to be obtained exactly. We
first give the following definition, which  is a generalized version
of Definition 2.2 (Khargonekar, Petersen and Zhou (1990).

{\bf Definition 12.} System (\ref{eq rob1}) (with $u=0$) is said to be quadratically stable,
if there exists  a positive definite matrix  $P>0$ and a constant $\alpha>0$,
such that the  differential generator of Lyapunov function $V(x)=x'Px$ satisfies
$$
{\cal L}[V(x)]=x'(PA+A'P+P\Delta A+(\Delta A)'P+C'PC)x\le -\alpha||x||^2
$$
for all $x\in {\cal R}^n$. (\ref{eq rob1}) is called quadratically
stabilizable, if there exists a state feedback $u(t)=Kx(t)$, such that the closed-loop
system is quadratically stable.

The following theorem is a necessary and sufficient condition for
quadratic stabilizability, and its proof needs two well known results as follows:

{\bf Lemma 2 (Xie, 1996).} For any real matrices $Y, H$ and $E$ of suitable dimensions,
$Y\in {\cal S}_n$, then for all $F$ with $F'F\le I$, we have
$$
Y+HFE+E'F'H'<0
$$
iff for some $\varepsilon>0$,
$$
Y+\varepsilon HH'+{\varepsilon}^{-1}E'E<0 .
$$

{\bf Lemma 3 (Schur's lemma).} For real matrices $N$, $M=M', R=R'>0$, the following two
conditions are equivalent:

1)\  $ M-NR^{-1}N'>0$ ,

2)\
$\left[
\begin{array}{cc}
M  &  N\\
N'  &  R
\end{array}
\right]>0.
$

{\bf Theorem 8.} System (\ref {eq rob1}) is quadratically
stabilizable iff there exist real matrices $Y$ and $X>0$, such that
\begin{equation}
\left[
\begin{array}{ccc}
AX+XA'+BY+Y'B'+CXC'+CY'D'+DYC' &  DY & XG'\\
Y'D' &  -X  & 0  \\
GX    &  0  &  -I
\end{array}
\right]<0,
\label {eq th8}
\end{equation}
especially,  $u(t)=Kx(t)=YX^{-1}x(t)$ is a quadratically stabilizing
control law.

{\bf Proof.}  By Definition 8, (\ref {eq rob1}) is quadratically stabilizable iff there exist
matrices $K$ and   $P>0$, such that
\begin{equation}
P(A+\Delta A+BK)+(A+\Delta A+BK)'P+(C+DK)'P(C+DK)<0.
\label {eq rob2}
\end{equation}
From Remark 1, $(A+\Delta A+BK,C+DK)$ is stable $\Leftrightarrow$ $(( A+\Delta A+BK)',(C+DK)')$
is stable, so by Lemma 1, (\ref{eq rob2}) is equivalent to
$$
P(A+\Delta A+BK)'+(A+\Delta A+BK)P+(C+DK)P(C+DK)'<0.
$$
Setting $Y_1=KP$, by Lemma 2, the above is equivalent to
\begin{equation}
\left[
\begin{array}{cc}
 P(A+\Delta A)'+(A+\Delta A)P+BY_1+Y'_1B'+CPC'+CY'_1D'+DY_1C'&  DY_1 \\
Y'_1D' &  -P
\end{array}
\right]<0.
\label {eq rob3}
\end{equation}
Take
$$
Z=\left[
\begin{array}{cc}
 PA'+AP+BY_1+Y'_1B'+CPC'+CY'_1D'+DY_1C' &  DY_1 \\
Y'_1D' &  -P
\end{array}
\right]<0,
$$
then (\ref{eq rob3}) can be written as
\begin{eqnarray}
Z&+&\left[
\begin{array}{cc}
\Delta A P+P(\Delta A)' &  0 \\
0 &  0
\end{array}
\right]=Z+\left[
\begin{array}{cc}
EFG P+PG'F'E' &  0 \\
0 &  0
\end{array}
\right] \nonumber\\
&=& Z+\left[
\begin{array}{c}
E \\
0
\end{array}
\right]
F\left[
\begin{array}{cc}
GP & 0
\end{array}
\right]+\left[
\begin{array}{c}
PG' \\
0
\end{array}
\right]F'\left[
\begin{array}{cc}
E' & 0
\end{array}
\right]<0.
 \label {eq rob5}
\end{eqnarray}
By Lemma 2, (\ref{eq rob5}) is equivalent to that for some $\varepsilon>0$,
\begin{eqnarray}
Z&+&\varepsilon\left[
\begin{array}{c}
E \\
0
\end{array}
\right]
\left[
\begin{array}{cc}
E' & 0
\end{array}
\right]+{\varepsilon}^{-1}\left[
\begin{array}{c}
PG' \\
0
\end{array}
\right]\left[
\begin{array}{cc}
GP & 0
\end{array}
\right]\nonumber \\
&=&\left[
\begin{array}{cc}
PA'+AP+BY_1+Y'_1B'+CPC'+CY'_1D'+DY_1C'\\
+\varepsilon EE'+{\varepsilon}^{-1} PG'GP &  DY_1 \\
Y'_1D' &  -P
\end{array}
\right]<0. \nonumber
\end{eqnarray}
The above divided by $\varepsilon$ with $X={\varepsilon}^{-1}P, Y={\varepsilon}^{-1} Y_1$,
results in
$$
\left[
\begin{array}{cc}
XA'+AX+BY+Y'B'+CXC'+CY'D'+DYC'\\
+EE'+ XG'GX &  DY  \\
Y'D' &  -X
\end{array}
\right]<0.
$$
Again, by Lemma 3,  the above is equivalent to
$$
\left[
\begin{array}{ccc}
AX+XA'+BY+Y'B'+CXC'+CY'D'+DYC' &  DY & XG'\\
Y'D' &  -X  & 0  \\
GX    &  0  &  -I
\end{array}
\right]<0.
$$
From our proof, $u(t)=Kx(t)=YX^{-1}x(t)$ is a quadratically
stabilizing contror law, the proof of Theorem 8 is completed.

 Since (\ref{eq th8}) is an LMI, by
some existing tools (Boyd, etal, 1994), one can easily test whether
it is empty or not, so Theorem 8 has practical value.

{\bf Remark 10.} By the same discussion as in Theorem 8, it is not
difficult to deal  with the quadratic stabilizability of
$$
dx=((A+\Delta A) x+(B+\Delta B)u)dt+((C+\Delta C)x+Du)dw,\ \ x(0)=x_0
$$
where uncertain matrices $\Delta A, \Delta B$ and  $\Delta C$ satisfy
$$
[\Delta A, \Delta B, \Delta C]=EF[G_1, G_2, G_3].
$$
Moreover, an analogous theorem to  Theorem 8 can be obtained.

\section{ Conclusion}

This paper studies the stabilizability (in mean square sense),
 weak stabilizability, exact observability of stochastic linear controlled
systems by the aid of spectrum technique.

Firstly,  some new concepts are introduced such as unremovable
spectrums, strong solutions of GAREs, etc.. Based on the spectrum
technique, we have obtained necessary and sufficient conditions for
both stabilizability and weak stabilizability of stochastic systems,
and have found  some new phenomena different from deterministic
systems. An important problem is that under what conditions with
$A,B,C$ and $D$, there always exists a feedback gain matrix  $K$,
such that for any given $\lambda _1, \lambda _2, \cdots, \lambda
_{n(n+1)/2} \in {\cal C}$, $\sigma({\cal L}^*_K)=\{\lambda _1,
\lambda _2, \cdots, \lambda _{n(n+1)/2}\}$?

Secondly, we should point out that Theorem 7 will have important
applications in the studies of GAREs, filtering, stochastic
stability.  further discussion will be included in our forthcoming
paper.

Finally, there are many topics in robust  quadratic stabilization of
stochastic systems needed to be studied, including time-varying (the
simplest case is $F=F(t)$) and nonlinear cases. All the above
problems are not only interesting but also important.

\end{document}